\documentclass[twocolumn]{autart}    

\usepackage[utf8x]{inputenc}
\usepackage{amsmath,amssymb}

\newtheorem{lemma}{Lemma}
\newtheorem{proposition}[lemma]{Proposition}
\newtheorem{theorem}{Theorem}
\newtheorem{corollary}[lemma]{Corollary}
\newtheorem{assumption}[lemma]{Assumption}

\newtheorem{definition}{Definition}
\newtheorem{remark}{Remark}
\begin{document}
\begin{frontmatter}

\title{Periodic excitations of  bilinear quantum systems}

\thanks{This paper was not presented at any IFAC
meeting.\\
 Tel: + 33 3 83 68 45 81,   Fax: + 33 3 83 68 45 34}

\author[Paestum]{Thomas Chambrion}\ead{Thomas.Chambrion@inria.fr},    

\address[Paestum]{IECN UMR 7502, Nancy University, CNRS, INRIA, BP 70239, 54506 Vand{\oe}uvre-l\`es-Nancy, France  and INRIA Nancy Grand Est, team CORIDA.}  %


\begin{abstract}                          
A well-known method of transferring the population of a quantum system from an eigenspace of the free Hamiltonian to another  is to use a periodic control law with an angular frequency equal to the difference of the eigenvalues.  For finite dimensional quantum systems, the classical theory of averaging provides a rigorous explanation of this experimentally validated result.   This paper extends this finite dimensional result, known as the  Rotating Wave Approximation,  to infinite dimensional systems and provides explicit convergence estimates.
\end{abstract}

\end{frontmatter}

\section{Introduction}

\subsection{Effective control of quantum systems}
The state of a quantum system evolving on a finite dimensional Riemannian manifold $\Omega$
is described by its \emph{wave function},
that is,
a point in the unit sphere of $L^2(\Omega, \mathbf{C})$.
In the absence of interaction with the environment and with a suitable choice of units, the time evolution of the
wave function is given by the Schr\"odinger equation
$$
\mathrm{i} \frac{\partial \psi}{\partial t}=-\frac{1}{2}\Delta \psi +V(x) \psi(x,t),
$$
where $\Delta$ is the Laplace-Beltrami operator on $\Omega$ (with suitable boundary conditions) and
$V:\Omega\rightarrow \mathbf{R}$ is a real function (usually called potential)
accounting for the physical properties of the system.  When subjected to an
excitation by an external electric field ({\it e.g.} a laser), the Schr\"odinger
equation reads
\begin{equation}\label{eq:blse}
\mathrm{i} \frac{\partial \psi}{\partial t}=-\frac{1}{2}\Delta \psi +V(x) \psi(x,t)
+u(t)
W(x) \psi(x,t),
\end{equation}
where $W:\Omega\rightarrow \mathbf{R}$ is a real function accounting for the
physical properties of the laser and $u$ is a real function of the time
accounting for the intensity of the laser.

A natural question, with many practical implications, is whether there exists a control $u$ that steers the quantum system  from a given initial wave function to a given target wave function (controllability issue) and, more important, how to build this control law (effective design of controls).

Considerable effort has been expended by different communities on studying the controllability of~(\ref{eq:blse}). We refer to  Nersessyan \cite{nersesyan}, Beauchard \& Mirrahimi \cite{beauchard-mirrahimi}, Mirrahimi \cite{mirrahimi-continuous}, Boscain \& Laurent \cite{camillo} and Boscain, Caponigro, Chambrion \& Sigalotti \cite{Schrod2} for a description of the known theoretical results concerning the existence of  controls steering a given source to a given target. As proved by Nersessyan \cite{nersesyan} and Mason \& Sigalotti  \cite{genericity-mario-paolo}, approximate controllability is a generic property for systems of the type (\ref{eq:blse}).

 A number of effective control algorithms have been obtained by various authors, see among many others Warren, Rabitz \& Dahleh \cite{warren},  Bl{\"{u}}mel, Fishman \& Smilansky \cite{Blumel}, Ohtsuki, Kono \& Fujiyama \cite {Ohtsuki} or Belhadj, Salomon \& Turinici \cite{MR2447854}. Most of the controls used in practice exhibit a remarkable pattern of periodic shape, with a frequency corresponding to the transitions of the quantum system (see for instance Salomon, Dion \& Turinici \cite{Salomon}).

\subsection{Averaging techniques in quantum mechanics}

The fact that a small amplitude periodic excitation with suitable frequency is sufficient, in general, to induce a transfer of quantum population from one energy level to another has prompted much attention.

The situation is now well understood for quantum systems with finitely many energy levels. These systems appear, for instance, as truncations of infinite dimensional systems. In this case, system \eqref{eq:blse} reads $\dot{x}=(A+uB)x$ where $A$ and $B$ are $N\times N$ skew-hermitian matrices. Let us briefly recall the method of the proof.

The mathematical concept of averaging of dynamical systems was introduced more than a century ago and has now developed into a well-established theory, see for instance the books of Guckenheimer \& Holmes \cite{Holmes}, Bullo \& Lewis \cite{Bullo} or Sanders, Verhulst \& Murdock \cite{Sanders}. It was observed that, for regular $F$ and small $\varepsilon$, the trajectories of the system $\dot{x}=\varepsilon F(x,t,\varepsilon)$
remain $\varepsilon$ close, for time of order $1/\varepsilon$, to the trajectories of the average system $\dot{x}=\widetilde{F}(x)$ where $\widetilde{F}(x)=\lim_{t\to \infty} 1/t \int_0^t F(x,t,0)$.

In quantum physics, this procedure has been known as the \emph{rotating wave approximation} for several decades, see Fox \& Eidson \cite{Fox} or Vandersypen \& Chuang \cite{NMR}. Assuming without loss of generality that $A$ is diagonal with eigenvalues $\mathrm{i}\lambda_1, \ldots, \mathrm{i}\lambda_N$, define, for every $n$ in $\mathbf{N}$, $y_n:t\mapsto e^{-tA} x_n$ where $x_n$ is the solution of $\dot{x}=(A+u^{\ast}(t)B/n)x$. The mapping $y_n$ is absolutely continuous and satisfies $\dot{y}_n=(u^{\ast}(t)/n)e^{-At}Be^{At}y_n$. The conclusion follows from standard averaging theory by computing the average matrix $\lim_{t\to \infty} 1/(nt) \int_0^t u^{\ast}(\tau)e^{-A\tau}Be^{A\tau} \mathrm{d}\tau$ whose all entries are zero but maybe entry $(j,k)$ if $u^{\ast}$ is $2\pi/|\lambda_j-\lambda_k|$-periodic.

By contrast with this finite dimensional result, the situation is much more intricate when the ambiant space has infinite dimension and the problem involves unbounded operators, which is precisely the case for the  Schr\"{o}dinger equation \eqref{eq:blse} when $\dim \Omega \geq 1$ and $A=\mathrm{i}\Delta/2$. The above averaging method fails because of serious regularity issues. For instance, the mapping $t\mapsto e^{tA}$ is no longer Lipschitz continuous.
  If moreover $B$ is unbounded, that is, if function $W$ in (\ref{eq:blse}) is not in $L^{\infty}(\Omega,\mathbf{R})$, then, for any given $t$, the mapping $x\mapsto e^{-tA}Be^{tA}x$ is not even continuous: this is what prevents the direct application of the averaging results in Banach spaces presented by Artstein \cite{Artstein}. For these reasons, most of the available averaging results for infinite dimensional quantum systems deal with constant controls only, as in the papers of Kummer \cite{Kummer} or Scherer \cite{Scherer}.

\subsection{Framework and notation}\label{SEC_notations}
 We first reformulate the problem \eqref{eq:blse} in a more abstract framework.  In a separable Hilbert space $H$ endowed with norm $\| \cdot \|$ and Hilbert product $\langle \cdot, \cdot \rangle$, we consider the following evolution problem:
\begin{equation}\label{EQ_main}
\frac{d \psi}{dt}=(A+u(t)B)\psi(t) \quad u(t) \in U
\end{equation}
where $(A,B,U)$ satisfies Assumption \ref{ASS_assumption_weak}.
\begin{assumption}\label{ASS_assumption_weak}
$(A,B,U)$ is a triple where $(A,B)$ is a pair of linear operators and $U$ is a subset of $\mathbf{R}$ such that
 \begin{enumerate}
 \item for every $n$ in $\mathbf{N}$, $U\subset n U$;
  \item $A$ is skew-adjoint with domain $D(A)$;
  \item there exists an Hilbert basis $(\phi_k)_{k \in \mathbf{N}}$ of $H$  and a family $(\mathrm{i}  \lambda_k)_{k \in \mathbf{N}}$ in $\mathrm{i} \mathbf{R}$ such that $A\phi_k=\mathrm{i}  \lambda_k \phi_k$ for every $k$ in $\mathbf{N}$;
  \item $B$ is skew-symmetric, possibly unbounded with domain $D(B)$;
  \item  for every $k$ in $\mathbf{N}$, $\phi_k$ belongs to $D(B)$;
\item for every $u$ in $U$, $A+uB$ is essentially skew-adjoint. \label{ASS_A+uB_skew-adjoint}
 \end{enumerate}
\end{assumption}

Assumption \ref{ASS_assumption_weak}.\ref{ASS_A+uB_skew-adjoint} ensures that, for every constant $u$ in $U$, $A+uB$ generates a group of unitary propagators. Hence, for every initial time $t_0$ in $\mathbf{R}$ and initial condition $\psi_0$ in $H$, for every piecewise constant control $u$ taking value in $U$, we can define the solution $t\mapsto \Upsilon^u(t,t_0) \psi_0$ of \eqref{EQ_main} taking value $\psi_0$ at time $t_0$. We simply note $\Upsilon^u(t,t_0) \psi_0= \Upsilon^u_t \psi_0$ when $t_0=0$.

\begin{remark}
To the best of the author's knowledge, in the general frame of Assumption \ref{ASS_assumption_weak}, no definition of solutions of \eqref{EQ_main} is available for  controls that are not piecewise constant.
With some extra regularity assumptions (for instance, when $B$ is bounded, see Section \ref{SEC_bounded}), one can define the solution of \eqref{EQ_main} for more general controls.
\end{remark}

\subsection{Main result}
\begin{definition} Let $(A,B,U)$ satisfy Assumption \ref{ASS_assumption_weak}.
A point $(j,k)$ of $\mathbf{N}^2$ is a \emph{non-degenerate transition} of $(A,B)$ if
(i) $j\neq k$, (ii) $\langle\phi_j,  B \phi_k\rangle \neq 0$ and (iii) for every $l,m$ in $\mathbf{N}$, $|\lambda_j-\lambda_k|=|\lambda_l-\lambda_m|$ implies $\{j,k\}=\{l,m\}$  or $\{j,k\} \cap \{l,m\}=\emptyset$ or  $\langle\phi_l,  B \phi_m\rangle = 0$.
\end{definition}

\begin{theorem}\label{PRO_main_result}
Let $(A,B,U)$ satisfy Assumption \ref{ASS_assumption_weak}, $(j,k)$ a non-degenerate transition of $(A,B)$ and $u^{\ast}:\mathbf{R}^+\rightarrow U$ be a piecewise constant function, periodic with period $T=\frac{2\pi}{|\lambda_j-\lambda_k|}$.
Assume that $\displaystyle{\int_0^T \!\!\! u^{\ast} (\tau) e^{  \mathrm{i}(\lambda_{l}-\lambda_{m})\tau} \mathrm{d}\tau = 0}$
for every $l,m$ such that $|\lambda_l-\lambda_m| \in (\mathbf{N}\setminus\{ 1 \}){|\lambda_j-\lambda_k|}$ and
$\langle \phi_l, B\phi_m\rangle \neq 0$ and $\{l,m\} \cap \{j,k\}\neq \emptyset$.

If $\displaystyle{\int_0^T \!\!\! u^{\ast} (\tau) e^{  \mathrm{i}(\lambda_{j}-\lambda_{k})\tau} \mathrm{d}\tau \neq 0}$,
  then there exists $T^{\ast}>0$ such that the sequence $\left (\left | \left \langle \phi_k, \Upsilon^{\frac{u^{\ast}}{n}}_{n T^{\ast}}(\phi_j)\right \rangle \right |\right )_{n \in \mathbf{N}} $ tends to $1$ as $n$ tends to infinity.
\end{theorem}

\subsection{Content of this paper}

This paper comprises three parts. The first one (Section \ref{SEC_finite_dimensional}) is concerned with a finite dimensional version of Theorem \ref{PRO_main_result}. As already mentioned, the convergence result is classical. A time reparametrization (Section \ref{SEC_time}) is introduced, which leads to explicit estimates better than those currently available in the literature.

The second part (Section \ref{SEC_general_case}) contains a general proof of Theorem \ref{PRO_main_result}, valid in the general framework of Assumption \ref{ASS_assumption_weak}. It provides the first available time estimates for the approximate controllability of general bilinear quantum systems when the free Hamiltonian admits a dense family of eigenvectors.

Finally, Section \ref{SEC_examples} presents an extension of Theorem \ref{PRO_main_result} to the case of the Morse quantum oscillator. The main feature of this system is that the spectrum of the free Hamiltonian has a continuous part.

\section{Finite dimensional estimates}\label{SEC_finite_dimensional}

\subsection{Finite dimensional framework}\label{SEC_finite_dimensional_framework}
Let $N$ be an integer and $A^{(N)}, B^{(N)}$ be two skew-adjoint matrices of order $N$.
Without lost of generality, we may assume that $A^{(N)}$ is diagonal with eigenvalues $\mathrm{i}\lambda_1,\mathrm{i}\lambda_2,\ldots,\mathrm{i}\lambda_N$. We denote with $(\phi_j)_{1\leq j\leq N}$ the canonical basis of $\mathbf{C}^N$, with $\langle \cdot,\cdot,\rangle$ the canonical Hermitian product of $\mathbf{C}^N$, with $\left (b_{jk}=\langle \phi_j, B^{(N)} \phi_k \rangle \right )_{1\leq j,k\leq N}$ the entries of $B^{(N)}$ and with $\pi_l^{(N)}:\mathbf{C}^N\rightarrow \mathbf{C}^N$ the projection $\pi_l^{(N)}:x\mapsto \sum_{m\leq l} \langle \phi_m,x\rangle \phi_m$ on the first $l$ components in $\mathbf{C}^N$.

For every piecewise constant function $u:\mathbf{R}\rightarrow \mathbf{R}$, we denote with $X_{(N)}^u(t,s)$ the propagator between times $s$ and $t$ of the system
\begin{equation} \label{EQ_Sigma_N}
 \dot{x}=(A^{(N)}+u(t)B^{(N)})x(t).
\end{equation}
In other words, $t\mapsto X_{(N)}^u(t,s) x_0$ is the unique solution of \eqref{EQ_Sigma_N} with initial condition $x(s)=x_0$. It is classical that, for every $t,s$ in $\mathbf{R}$, for every $x_0$ in $\mathbf{C}^N$, the mapping $u\mapsto  X_{(N)}^u(t,s)x_0$ admits a
Lipschitz-continuous continuation to $L^1_{loc}(\mathbf{R},\mathbf{R})$.

The aim of this section is to prove the following result.
\begin{proposition}\label{PRO_estimates_dim_finie}
Let $u^{\ast}:\mathbf{R}^+\rightarrow \mathbf{R}$ be a locally integrable function.
Assume that $\lambda_1 \neq \lambda_2$ and that, for every $l,m\leq N$, $|\lambda_l-\lambda_m|=|\lambda_2-\lambda_1|$ implies $\{l,m\}=\{1,2\}$ or $b_{lm}=0$ or $\{l,m\}\cap \{1,2\}=\emptyset$.

 Assume that $u^{\ast}$ is periodic with period $T=\frac{2\pi}{|\lambda_2-\lambda_1|}$  and that $\displaystyle{\int_0^T \!\!\! u^{\ast} (\tau) e^{  \mathrm{i}(\lambda_{l}-\lambda_{m})\tau} \mathrm{d}\tau = 0}$
for every $\{l,m\}$ such that $\{l,m\}\cap \{1,2\}\neq \emptyset$ and ${\lambda_l-\lambda_m}\in (\mathbf{Z}\setminus\{\pm1\})(\lambda_1-\lambda_2)$ and
$b_{lm} \neq 0$.

If $\displaystyle{\int_0^T \!\!\! u^{\ast} (\tau) e^{  \mathrm{i}(\lambda_{2}-\lambda_{1})\tau} \mathrm{d}\tau \neq 0}$,
  then, for every $n$ in $\mathbf{N}$, there exists $T_n^{\ast}$ in $(nT^{\ast}-T,nT^{\ast}+T)$ such that
\begin{equation}\label{EQ_estimation_convergence}
\frac{1-|\langle \phi_2,X^{\frac{u}{n}}_{(N)}(T^{\ast}_n,0) \phi_1 \rangle |}{ (1+ 2 K\| B^{(N)}\|)I}\leq
 \frac{ (1+C)\|\pi_2^ {(N)} B^{(N)} \|}{n},
\end{equation}
with 
$$T^{\ast}=  \frac{\pi  T}{2 |b_{1,2}|  \left |\int_0^T \!\! \!u^{\ast}(\tau)e^{\mathrm{i}(\lambda_{1}-\lambda_{2}) \tau}  \mathrm{d}\tau \right |},\quad I=\int_0^T \! \!|u^{\ast}(\tau)|\mathrm{d}\tau,$$
$$ K=\frac{IT^{\ast}}{T}, C=\sup_{(j,k)\in \Lambda} \left | \frac{\int_0^T u^{\ast}(\tau) e^{\mathrm{i} (\lambda_j-\lambda_k)\tau}\mathrm{d}\tau}{\sin \left ( \pi\frac{|\lambda_j-\lambda_k|}{|\lambda_2-\lambda_1|} \right )} \right |,$$
where $\Lambda$ is the set of all pairs $(j,k)$ in  $\{1,\ldots,N\}^2$  such that $b_{jk} \neq 0$ and $\{j,k\}\cap\{1,2\} \neq \emptyset$ and $ |\lambda_j-\lambda_k|\notin \mathbf{Z}|\lambda_2-\lambda_1|$.
\end{proposition}
\begin{corollary}\label{COR_convergence}
With the notations of Proposition \ref{PRO_estimates_dim_finie}, $|\langle \phi_2,X^{\frac{u}{n}}_{(N)}(n T^{\ast},0) \phi_1 \rangle |$ tends to 1 as $n$ tends to infinity.
\end{corollary}
The convergence result (Corollary \ref{COR_convergence}) is classical. The novel element here is that the bilinear structure allows us to give explicit estimates (\ref{EQ_estimation_convergence}) for the convergence rate.

Most of the rest of this section is devoted to the proof of Proposition \ref{PRO_estimates_dim_finie}. A technical time reparametrization (Section \ref{SEC_time}) is introduced. The main ingredient of the proof is an averaging procedure (Section \ref{SEC_averaging_procedure}), that is performed first on piecewise constant controls and then extended to irregular controls. Finally, we introduce the notion of \emph{efficiency} in Section \ref{SEC_Efficiency}.

\subsection{Time reparametrization}\label{SEC_time}

We note $PC$ the set of the piecewise constant functions for which  there exist two  sequences $(u_j)_{1\leq j \leq p}$ and $(t_j)_{1\leq j \leq p}$ with value in $(0,+\infty)$ such that
$$
u=\sum_{1\leq j \leq p+1} u_j \chi_{[\tau_j,\tau_j+t_j)},
$$
where $\chi$ is the characteristic function and the sequence $(\tau_j)_{1\leq j\leq p+1}$ is defined by induction: $\tau_1=0$, $\tau_{j+1}=\tau_j +t_j$.
An element $u$ of $PC$ will be denoted $(u_j,t_j)_{1\leq j \leq p}$.

The involutary mapping $ {\mathcal P}:(u_j,t_j)_{1\leq j \leq p}\in PC \mapsto   \left ( {1}/{u_j}, u_j t_j \right )_{1\leq j \leq p}  \in PC$
can be used as a time reparametrization to replace the system (\ref{EQ_Sigma_N})
 by the control system in $\mathbf{C}^N$
\begin{equation}\label{EQ_main_check}
 \frac{\mathrm{d} x}{\mathrm{d} t}=(u(t) A^{(N)} +B^{(N)})x(t),
\end{equation}
whose propagator between time $s$ and $t$  will be denoted with $\check{X}_{(N)}^u(t,s)$.

\begin{proposition}\label{PRO_reparam_temps}
 For every $x_0$ in $\mathbf{C}^N$ and every $u$ in $PC$, $\check{X}_{(N)}^{\mathcal{P}u}({\int_0^t u(\tau)\mathrm{d}\tau},0)x_0=X_{(N)}^u(t,0)x_0$.
\end{proposition}
\textbf{Proof:}
  This follows from the equality $e^{t(A^{(N)}+uB^{(N)})}=e^{tu \left (\frac{1}{u} A^{(N)}+B^{(N)} \right )}$, valid on every interval where $u$ is constant.
\qed

\subsection{Averaging procedure}\label{SEC_averaging_procedure}

Let $u^{\ast}$ be a non vanishing piecewise constant function, $T$-periodic, as in the hypotheses of Proposition~\ref{PRO_estimates_dim_finie}.
For every $n$ in $\mathbf{N}$, define the non-vanishing $T$-periodic function $u_n= u^{\ast}/n$. For every $n$ in $\mathbf{N}$, $\int_0^{nT} |u_n(\tau)| \mathrm{d} \tau=\int_0^T |u^{\ast}(\tau)| \mathrm{d}\tau$.

We perform now the time reparametrization set out in Section \ref{SEC_time}.
The function $t\mapsto \int_0^t |u_n(\tau)| \mathrm{d}\tau$ is non-decreasing. We denote with $v_n$ its reciprocal function and define also $I=\int_0^T \!\! |u^{\ast}(\tau)| \mathrm{d}\tau$. For every $t$ in $\mathbf{R}^+$, $v_n \left (t+{I}/{n} \right )=v_n(t)+T$ and
${\mathcal P} |u_n|$ is the ${I}/{n}$ periodic  derivative (defined almost everywhere) of $v_n$.

Fix $x_0$ in $\mathbf{C}^N$ and note with $t\mapsto x_n(t)$  the solution of
\begin{equation}\label{EQ_dyn_x}
\dot{x}=(A^{(N)}+u_n(t)B^{(N)})x(t)
\end{equation}
 with initial condition $x(0)=x_0$.

The set
$[0,T]$ can be written as a finite union of  disjoint intervals
$$[0,T]=J_1^{+}\cup J_1^{-} \cup \ldots \cup J_p^{+}  \cup J_p^{-},$$
such that $u^{\ast}$ takes positive (respectively, negative) values on $J^{+}=\cup_{l=1}^p J_l^{+}$ (respectively,
 $J^{-}=\cup_{l=1}^p J_l^{-}$). Defining the sets  $G^{+}_n=  v_n ^{[-1]}(J^{+})$
$=\{l \in \mathbf{R}^{+}| \exists s \in J^{+} , \int_0^s |u_n(\tau)|\mathrm{d}\tau= l\}$ and
$G^{-}_n= v_n^{[-1]}(J^{-})$ $=\{l \in \mathbf{R}^{+}| \exists s \in J^{-} , \int_0^s |u_n(\tau)|\mathrm{d}\tau= l\}$, we obtain the dynamics of $y_n= x_n \circ v_n$, valid for almost every $t$:
\begin{equation}\label{EQ_dyn_yn}
\frac{\mathrm{d}y_n}{\mathrm{d}t}=
\left \{ \begin{array}{ll}({\mathcal P}(|u_n|)(t) A^{(N)}+B^{(N)})y_n(t) & \mbox{if } t \in G_n^{+}\\
			  ({\mathcal P}(|u_n|)(t) A^{(N)}-B^{(N)})y_n(t)	& \mbox{if } t \in G_n^{-}
         \end{array} \right.
 \end{equation}
Finally,  for every $t$, we define $z_n(t)=e^{-v_n(t)A}y_n(t)$
(note that, for every $t$, for every $l$ in $\mathbf{N}$, $|\langle \phi_l,z_n(t)\rangle|=
 |\langle \phi_l, y_n(t) \rangle |$), and the time varying $N \times N$ matrix $M_n$ 
$$M_n:t\mapsto \mathrm{sg}(u_n\circ v_n)  e^{-v_n(t)A^{(N)}}B^{(N)}e^{v_n(t)A^{(N)}}.$$
From (\ref{EQ_dyn_yn}), we deduce the dynamics of $z_n$, valid for almost every $t$ in $G^{+}_n\cup G^{-}_n$:
\begin{equation}\label{EQ_dyn_zt}
\frac{\mathrm{d} z_n}{\mathrm{d} t}\!=M_n(t) z_n(t).
\end{equation}
We note $Z^n_t$ the propagator associated with \eqref{EQ_dyn_zt}.
Note that, for every $k$, the mapping  $t\mapsto \langle \phi_k,z_n\rangle$ is Lipschitz continuous with Lipschitz constant $\|B \phi_k\|$.

Let $M^{\dag}$ be the constant $N\times N$ matrix whose entries,  for $1 \leq j,k \leq N$,  are defined by
$$m^{\dag}_{j,k}= \frac{b_{j,k}}{I} \int_0^I \exp \left (\mathrm{i} (\lambda_{j}-\lambda_{k}) v^{\ast}(\tau) \right )\mathrm{d}\tau $$ if $T (\lambda_{j}-\lambda_{k}) \in 2 \pi \mathbf{Z}$ and $m^{\dag}_{j,k}=0$ if $T (\lambda_{j}-\lambda_{k}) \notin 2 \pi \mathbf{Z}$. Define the $N\times N$ matrix $H_n(t)$, with entries $(h_{jk}^n(t))_{1\leq j,k\leq N}$, by $H_n(t)=M_n(t)-M^{\dag}$.
\begin{lemma}\label{LEM_lemme_inegalites_hij}
 For every $t$ in $\mathbf{R}^+$, for every $n$ in $\mathbf{N}$, for every $j,k\leq N$ such that $|\lambda_j-\lambda_k| \notin \mathbf{Z}|\lambda_2-\lambda_1|$,
 \begin{equation}\label{EQ_inegalite_non_resonant}
 \left |\int_0^t h_{jk}^n(\tau)\mathrm{d}\tau \right | \leq
 \frac{|b_{jk}|}{n}\left (\left |\frac{\int_0^T u^{\ast}(\tau) e^{\mathrm{i}\tau(\lambda_j-\lambda_k)}\mathrm{d}\tau}{\sin \left ( \pi \frac{\lambda_j-\lambda_k}{\lambda_2-\lambda_1} \right )} \right | + I\right ).
\end{equation}
 For every $t$ in $\mathbf{R}^+$, for every $n$ in $\mathbf{N}$, for every $j,k\leq N$ such that $\{j,k\}\neq \{1,2\}$ and $|\lambda_j-\lambda_k|/|\lambda_2-\lambda_1| \in \mathbf{Z}$,
 \begin{equation}\label{EQ_inegalite_resonant}
 \left |\int_0^t h_{jk}^n(\tau)\mathrm{d}\tau  \right | \leq
 \frac{2|b_{jk}|I}{n}.
 \end{equation}
\end{lemma}
\textbf{Proof:}
For every $t$, define the integer $s=\lfloor \frac{tn}{I} \rfloor$.  For every $j$, $k$, $n$, $t$, we have
$$
\left |\int_{s I/n}^t \!\!\!\!\! \mathrm{sg}(u_n\circ v_n)e^{ \mathrm{i}  (\lambda_j-\lambda_k) v_n(\tau)}\mathrm{d}\tau\right |\leq \frac{I}{n}.
$$
 For every $j,k$ such that $T (\lambda_{j}-\lambda_{k}) \notin 2 \pi \mathbf{Z}$,
\begin{eqnarray*}
\lefteqn{\left |\int_0^{s I/n} \!\!\!\!\!\!\! \!\!\mathrm{sg}(u_n\circ v_n)e^{ \mathrm{i}(\lambda_{j}-\lambda_{k}) v_n(\tau)}\mathrm{d}\tau \right |}\\
&\leq & \left |\sum_{m=1}^s \int_0^ {I/n} \!\!\!\!\!\!\! \mathrm{sg}(u_n\circ v_n)e^{\mathrm{i}(\lambda_{j}-\lambda_{k}) (v_n(\tau)+mT)} \mathrm{d}\tau \right | \\
&\leq & \left | \int_0^ {I/n} \!\!\!\!\!\!\! \mathrm{sg}(u_n\circ v_n)e^{ \mathrm{i}(\lambda_{j}-\lambda_{k}) v_n(\tau)}\sum_{m=1}^s e^{\mathrm{i} (\lambda_{j}-\lambda_{k})  mT }\mathrm{d}\tau \right | \\
&\leq & \frac{1}{n}\left ( \frac{2\left |\int_0^T \!\!\!u^{\ast}(\tau) e^{\mathrm{i}(\lambda_{j}-\lambda_{k}) \tau} \mathrm{d}\tau\right |}{|1-\exp(\mathrm{i}T(\lambda_{j}-\lambda_{k}))|}  \right ).
\end{eqnarray*}

If $T({\lambda_{j}-\lambda_{k}}) \in 2 \pi \mathbf{Z}$, then
\begin{eqnarray*}
\lefteqn{\sum_{m=1}^s \int_0^ {I/n}\!\!\!\! \!\!\!\mathrm{sg}(u_n\!\!\circ v_n) e^{\mathrm{i} (\lambda_{j}-\lambda_{k}) (v_n(\tau)+mT) }\mathrm{d}\tau}\\  &\quad \quad \quad=&s
\int_0^ {I/n}\!\!\!\!\!\! \mathrm{sg}(u_n\!\!\circ v_n)e^{\mathrm{i} (\lambda_{j}-\lambda_{k}) v_n(\tau)} \mathrm{d}\tau\\
&\quad \quad \quad=&\frac{s}{n} \int_0^I \!\!\!\mathrm{sg}(u^{\ast}\!\!\circ v^{\ast}) e^{\mathrm{i} (\lambda_{j}-\lambda_{k}) v^{\ast}(\tau)} \mathrm{d}\tau.
\end{eqnarray*}
Hence
\begin{eqnarray*}
\lefteqn{\left |\int_0^{sI/n}\!\!\!\!\!\!\mathrm{sg}(u_n\circ v_n) e^{\mathrm{i} (\lambda_{j}-\lambda_{k}) v_n(\tau)}\mathrm{d}\tau \right.}\\
&&\quad \quad \quad - \left.\frac{t}{I} \int_0^T \!\!\!u^{\ast}(\tau) e^{\mathrm{i}(\lambda_{j}-\lambda_{k}) \tau} \mathrm{d}\tau \right | \leq  \frac{I}{n}.\qed
\end{eqnarray*}

Recall that, for every $n$ in $\mathbf{N}$ and every $t\geq 0$,
$$\dot{z}_n=M^{\dag}z_n +H_n(t) z_n,$$
or, using the variation of the constant formula,
$$z_n(t)=e^{tM^{\dag}}z_n(0)+ \int_0^t e^{(t-s)M^{\dag}} H_n(s) z_n(s) \mathrm{d}s.$$
Integrating by part, we get for every $t\geq 0$,
\begin{eqnarray}
z_n(t)&=&e^{tM^{\dag}}z_n(0)+\left \lbrack e^{(t-s)M^{\dag}} \left (\int_0^sH_n(\tau)\mathrm{d}\tau \right) z_n(s) \right \rbrack_0^t \nonumber\\
&& \quad +\int_0^t M^{\dag} e^{(t-s)M^{\dag}} \left (\int_0^sH_n(\tau)\mathrm{d}\tau \right) z_n(s) \mathrm{d}s \nonumber \\&& \quad - \int_0^t  e^{(t-s)M^{\dag}} \left (\int_0^sH_n(\tau)\mathrm{d}\tau \right) \dot{z}_n(s) \mathrm{d}s. \label{EQ_majoration_integration_par_partie}
\end{eqnarray}
The equality $X^{u_n}_{(N)}(t,0)=e^{tA^{(N)}}\circ Z^{u_n}_{(N)}(v_n^{[-1]}(t),0)$ gives  estimates for the convergence of the propagators $X^{u_n}$: for every $n$ in $\mathbf{N}$, for every $t\leq v_n^{[-1]}(K)$, for every $x_0$ in $\mathbf{C}^N$, ($I$, $K$ and $C$ are defined as in Proposition \ref{PRO_estimates_dim_finie})
\begin{eqnarray}
 \lefteqn{\frac{\|X^{u_n}_{(N)}(t,0)x_0-e^{tA^{(N)}}e^{v_n^{[-1]}(t)M^{\dag}}x_0\|}{I(C+1)\|M^{\dag}\|} }\nonumber\\ &&\leq
\frac{1+ K (\|M^{\dag}\|+ \sup_{s\leq t}\|B^{(N)}X^{u_n}_{(N)}(s,0)x_0\|)}{n}
\label{EQ_estimation_conv_uniform_propaga_dim_finie_weakly_coupled},
\end{eqnarray}
or
\begin{eqnarray}\label{EQ_major_dist_propagateur}
\frac{\|X^{u_n}_{(N)}(t,0)-e^{tA^{(N)}}e^{v_n^{[-1]}(t) M^{\dag}}\|}{I(C+1)  \|B^{(N)}\|} \leq
\frac{ 1 + 2K \|B^{(N)}\|}{n}\label{EQ_estimation_conv_uniform_propaga_dim_finie}.
\end{eqnarray}

Projecting (\ref{EQ_majoration_integration_par_partie}) on the two first components of $\mathbf{C}^N$, we get
\begin{eqnarray}\label{EQ_majoration_tout_temps_reparam_dim_finie}
\lefteqn{\left \| \pi_2^{(N)} z_n(t)- \pi_2^{(N)} e^{tM^{\dag}}z_n(0) \right \|}\nonumber \\
&\leq &
\left \| \pi_2^{(N)}\! \int_0^t \!\!H_n(\tau)\mathrm{d}\tau \right \| +t \left \| \pi_2^{(N)}\! M^{\dag} \right \|  \left \| \pi_2^{(N)}\! \int_0^t \!\!H_n(\tau)\mathrm{d}\tau \right \| \nonumber \\&&\quad + t\left \| \pi_2^{(N)} \int_0^t\!\! H_n(\tau)\mathrm{d}\tau \right \|  \| B^{(N)} \|.
\end{eqnarray}
Note that $|\langle \phi_2,e^{K M^{\dag}}\phi_1\rangle |=1$. Equation (\ref{EQ_estimation_convergence}) follows from \eqref{EQ_majoration_tout_temps_reparam_dim_finie} with $t=K$ and
$$ T^{\ast}_n= v_n (K)=v^{\ast} \left (\frac{n \pi I }{2 |b_{1,2}| \left |\int_0^T \! u^{\ast}(\tau) e^{\mathrm{i} (\lambda_{2}-\lambda_{1}) \tau} \mathrm{d}\tau \right |} \right ).$$

Knowing that $v^{\ast}$ is non-decreasing and $v^{\ast}(l I)= lT$ for every $l $ in $\mathbf{N}$, we deduce that, for every $n$ in $\mathbf{N}$, $nT^{\ast}-T\leq T_n^{\ast}\leq nT^{\ast}+T$, with $T^{\ast}$ defined as in the statement of Proposition~\ref{PRO_estimates_dim_finie}. 

Note finally that, for every $s,t$ in $\mathbf{R}$ such that $s\leq t$,
$$ \big | |\langle \phi_2, x_n(t) \rangle | -|\langle \phi_2, x_n(s) \rangle | \big |\leq \frac{\|B \phi_2\|}{n}\int_s^t|u(\tau)\mathrm{d}\tau.$$
 This completes the proofs of Proposition \ref{PRO_estimates_dim_finie} and Corollary \ref{COR_convergence} in the case where $u^{\ast}$ is a non vanishing piecewise constant function.

If $u^{\ast}$ is a locally integrable $T$-periodic function, let $(u^{\ast,l})_{l \in \mathbf{N}}$ be a sequence of non-vanishing piecewise constant $T$-periodic functions converging to $u^{\ast}$ in the distribution sense with $\|u^{\ast,l}\|_{L^1([0,T])}\leq |u|([0,T])$ for every $l$. We define, for every $n$, $u_{n,l}=\frac{u^{\ast,l}}{n}$, $v_{n,l}:t\mapsto\int_0^t {\mathcal P}|u_{n,l}|(\tau) \mathrm{d}\tau$ and  $M_{n,l}(t)=\mathrm{sg}(u_{n,l}\circ v_{n,l}(t))e^{-v_{n,l}(t) A}B^{(N)}e^{v_{n,l}(t) A}$. For every $t$, the matrix $ \int_0^tM_{l,l}(\tau) \mathrm{d}\tau$ tends to $t M^{\dag}$, uniformly with respect to $t$ in a compact set, as  $l$ tends to infinity. Hence the solutions of $\dot{x}=M_{l,l}(t)x$ tend to the solutions of $\dot{x}=M^{\dag} x$, uniformly with respect to the time in a compact interval, as $l$ tends to infinity. That concludes the proof of Proposition~\ref{PRO_estimates_dim_finie}.

\subsection{Efficiency of the transfer}\label{SEC_Efficiency}

Continuing with the notation of the last paragraph, for every non identically zero, $\frac{2\pi}{|\lambda_j-\lambda_k|}$-periodic function $u^{\ast}$, we define the efficiency of $u^{\ast}$ with respect to the transition $(j,k)$ as the real quantity:
\begin{eqnarray*}
E^{(j,k)}(u^{\ast})&=&
\frac{\left |\int_0^{\frac{2\pi}{|\lambda_j-\lambda_k|}}\! u^{\ast}(\tau) e^{\mathrm{i}(\lambda_j-\lambda_k) \tau} \mathrm{d}{\tau} \right |}{\int_0^{\frac{2\pi}{|\lambda_j-\lambda_k|}}\! |u^{\ast}(\tau) | \mathrm{d}{\tau} } .
\end{eqnarray*}
For every $u$, $0\leq E^{(j,k)}(u) \leq 1$.
For every $\{j, k\}$, $\sup_u  E^{(j,k)}(u)=1$. The supremum is reached with a periodic sum of Dirac masses.

An intuitive explanation of the efficiency could be the following: asymptotically, the $L^1$ norm of the control  needed to induce the transition between levels $j$ and $k$ using periodic controls of the form $u_n$ is equal to
$\pi/(2|b_{j,k}| E^{(j,k)}(u^{\ast}))$.

The system \eqref{EQ_main} being given, the design of an effective control law fulfilling the hypotheses of Theorem~\ref{PRO_main_result}  is an important
 practical issue.
 To generate  a transfer from level $j$ to level $k$, one should choose a control $u$ such that $E^{(j,k)}(u)$ be as large as possible and $E^{(l_1,l_2)}(u)$ be zero (or arbitrarily close to zero) for every $l_1,l_2$ such that $\lambda_{l_1}-\lambda_{l_2} \in (\lambda_j-\lambda_k) \mathbf{Z}$. The algorithm we have described in \cite{Schrod2} allows us to build $u$ such that $E^{(j,k)}(u)>0.43 $, with $E^{(l_1,l_2)}(u)$ arbitrarily small for every finite number of pairs $\{l_1,l_2\}$ satisfying $\{l_1,l_2\}\neq \{j,k\}$ and  $|\lambda_{l_1}-\lambda_{l_2}|\neq |\lambda_{j}-\lambda_{k}|$. Some other examples, including also examples of ineffective controls with zero  efficiency, are studied in \cite{FEPS_ACC}.

\section{Infinite dimensional estimates}\label{SEC_general_case}

We come back to the general case of Assumption \ref{ASS_assumption_weak}. To ensure that the system (\ref{EQ_main}) is well-posed, we consider only piecewise constant functions $u^{\ast}$. The method of the proof is directly inspired by \cite{Schrod2}: the original infinite dimensional system \eqref{EQ_main} is approached by a suitable Galerkin approximation, which allows us to apply the finite dimensional results of Section \ref{SEC_finite_dimensional}.

\subsection{Galerkin approximation}\label{SEC_Galerkin_approxim}
Let $(A,B,U)$ satisfy Assumption \ref{ASS_assumption_weak}. We define $b_{jk}=\langle \phi_j,B\phi_k\rangle$ for $j,k$ in $\mathbf{N}$. For every $N$ in $\mathbf{N}$, we define $\pi_N: \psi\in H\mapsto \sum_{j\leq N}\langle \phi_j,\psi\rangle \phi_j \in H$ and the compressions $A^{(N)}=\pi_N \circ A \circ \pi_N$ and $B^{(N)}=\pi_N\circ B \circ \pi_N$. Note that $A^{(N)}$ and $B^{(N)}$ are finite rank operators defined in an infinite dimensional space.  With an obvious abuse of notation, we extend the finite dimensional propagator  $X_{(N)}^u$  defined in Section \ref{SEC_finite_dimensional_framework} to the  infinite dimensional space $H$.

Let $u^{\ast}:\mathbf{R}\rightarrow U$, $j,k$ be given as in the hypotheses of Theorem \ref{PRO_main_result}. Up to a reordering, we may assume $j=1$ and $k=2$ without loss of generality.

Define $$K=\frac{\pi}{2 |b_{12}|} \frac{1}{Eff^{(1,2)}(u^{\ast})}.$$
Fix $\varepsilon>0$. Since $\phi_1$ and $\phi_2$ belong to the domain of $B$, the sequences $(b_{1,l})_{l\in \mathbf{N}}$ and
$(b_{2,l})_{l\in \mathbf{N}}$ are  in $\ell^2$. Hence, there exists $N$ in $\mathbf{N}$ such that $\|\pi_2 B(1-\pi_N)\|=\|(1-\pi_N) B\pi_2 \|<\varepsilon/K.$ 

Consider system (\ref{EQ_main}) with control $u_n:=u^{\ast}/n$ in projection   on $\mathrm{span}(\phi_1,\ldots,\phi_N)$:
\begin{eqnarray}
 \pi_N \frac{\mathrm{d}}{\mathrm{d}t} \Upsilon^{u_n}_t \phi_1 &=&
(A^{(N)} + u_n(t) B^{(N)}) \Upsilon^{u_n}_t \phi_1 \nonumber \\&&+u_n(t) \pi_N B (1-\pi_N) \Upsilon^{u_n}_t \phi_1.
\end{eqnarray}
From the variation of the constant, we get
\begin{eqnarray}
\lefteqn{\pi_N \Upsilon^{u_n}_t \phi_1= X_{(N)}^{u_n}(t,0) \phi_1} \nonumber\\
&&+ \int_0^t\!\!\!u_n(s) X^{u_n}_{(N)}(t,s) \pi_N B (1-\pi_N) \Upsilon^{u_n}_t \phi_1 \mathrm{d}s.\label{EQ_proj_N}
\end{eqnarray}
Project (\ref{EQ_proj_N}) on $\mathrm{span}(\phi_1,\phi_2)$, and recall that $\pi_N \pi_2=\pi_2 \pi_N=\pi_2$ for $N\geq 2$:
\begin{eqnarray}
 \lefteqn{\pi_2 \Upsilon^{u_n}_t \phi_1= \pi_2 X_{(N)}^{u_n}(t,0) \phi_1} \nonumber \\&&+ \int_0^t \!\!\!u_n(s) \pi_ 2 X_{(N)}^{u_n}(t,s) \pi_N B (1-\pi_N) \Upsilon^{u_n}_t \phi_1 \mathrm{d}s.\label{EQ_pre_bracket}
\end{eqnarray}
Define, for every $t,s$ in $\mathbf{R}$, the bounded linear mapping $[\pi_2, X^{(N)}(t,s)]:=\pi_2\circ X_{(N)}^{u_n}(t,s)-X_{(N)}^{u_n}(t,s) \circ \pi_2$. Equation (\ref{EQ_pre_bracket}) reads
\begin{eqnarray}
\lefteqn{\pi_2 \Upsilon^{u_n}_t \phi_1- \pi_2 X_{(N)}^{u_n}(t,0) \phi_1 =}\nonumber \\
&&
-\int_0^t\!\!\! u_n(s)  X_{(N)}^{u_n}(t,s) \pi_2 B (1-\pi_N) \Upsilon^{u_n}_t \phi_1 \mathrm{d}s  \nonumber \\
&& +
\int_0^t\!\!\! u_n(s) [\pi_2, X_{(N)}^{u_n}(t,s)]  \pi_N B (1-\pi_N) \Upsilon^{u_n}_t \phi_1 \mathrm{d}s.\label{EQ_post_bracket}
\end{eqnarray}

\subsection{Estimates of commutators}
Extend the definition of $M^{\dag}$ of Section \ref{SEC_averaging_procedure} by $M^{\dag}\phi_j=0$ for $j>N$ and
define the linear operator $E^n_{N}(t):=X_{(N)}^{u_n}(t,0)-e^{v^{[-1]}(t)M^{\dag}}\!\!$.
 Since the commutator $[\pi_2,M^{\dag}]=\pi_2 M^{\dag}-M^{\dag}\pi_2$ vanishes, we have, for every $t$ in $\mathbf{R}$,
\begin{eqnarray*}
\|[\pi_2, X_{(N)}^{u_n}(t,0)]\|&=&\|[\pi_2, e^{v^{[-1]}(t)M^{\dag}} +E^n_{(N)}(t)]\|\\
&=&\|[\pi_2,E^n_{(N)}(t)]\|\leq 2 \|E^n_{(N)}(t)\|.
\end{eqnarray*}
Note also that, for every $t$ in $\mathbf{R}$,
\begin{eqnarray*}
\|[\pi_2, X_{(N)}^{u_n}(0,t)]\|&=&
\|X_{(N)}^{u_n}\!(0,t) [X_{(N)}^{u_n}\!(t,0), \pi_2] X_{(N)}^{u_n}\!(0,t)    \|\\
&\leq& 2 \|E^n_{(N)}(t)\|.
\end{eqnarray*}
For every $s,t$ in $\mathbf{R}$,
\begin{eqnarray*}
 \lefteqn{[\pi_2, X_{(N)}^{u_n}(t,s)]}\\
&=&\pi_ 2X_{(N)}^{u_n}(t,0)X_{(N)}^{u_n}(0,s)-X_{(N)}^{u_n}(t,0)X_{(N)}^{u_n}(0,s) \pi_2\\
&=&X_{(N)}^{u_n}(t,0)[\pi_ 2,\! X_{(N)}^{u_n}(0,s)]+ [\pi_2, X_{(N)}^{u_n} (t,0)]X_{(N)}^{u_n}(0,s).
\end{eqnarray*}
Finally, we get, for every $(s,t)$ in $\mathbf{R}$, for every $n,N$ in $\mathbf{N}$.
\begin{equation}\label{EQ_maj_commutateurs}
 \left \|\left [ \pi_2, X_{(N)}^{u_n}(t,s) \right \rbrack \right \|\leq  4 \|E_{(N)}^n(t)\|.
\end{equation}

\subsection{Proof of Theorem \ref{PRO_main_result}}
From (\ref{EQ_post_bracket}) and (\ref{EQ_maj_commutateurs}), since $\|\pi_2B(1-\pi_N)\|<\varepsilon/K$,
\begin{eqnarray}
 \lefteqn{|\langle \phi_2,\Upsilon^{\frac{u^{\ast}}{n}}(t)\phi_1 \rangle -\langle \phi_2,X^{\frac{u^{\ast}}{n}}_{(N)}(t,0)\phi_1 \rangle|} \nonumber\\ &\quad \quad \quad \quad \quad \quad \leq&
\varepsilon  + 4 \|E_{(N)}^n(t)\| K \|\pi_N B(1-\pi_N)\|.\label{EQ_majoration_generale}
\end{eqnarray}
From \eqref{EQ_major_dist_propagateur}, $\sup_{t\leq v_n(K)} \|E^n_{(N)}(t)\|$ tends to zero as $n$ tends to infinity.
For $n$ large enough, $$\|E_{(N)}^n(nT^{\ast})\|\leq  \frac{\varepsilon}{4K \|\pi_N B(1-\pi_N)\|},$$ and $1-|\langle \phi_2,\Upsilon^{\frac{u^{\ast}}{n}}(nT^{\ast})\phi_1 \rangle |\leq 2 \varepsilon$. This completes the proof of Theorem \ref{PRO_main_result}.

\subsection{Bounded coupling}\label{SEC_bounded}

In this Section, we consider the special case where $B$ is bounded. The propagator $u\mapsto \Upsilon^u$ is defined on the set of locally integrable functions.

\begin{proposition}\label{PRO_cosinus}
Let $(A,B,\mathbf{R})$ satisfy Assumption \ref{ASS_assumption_weak}.
Assume that $B$ is bounded  and that $(1,2)$ is a non degenerate transition of $(A,B)$.
Define $T=\frac{2\pi}{|\lambda_2-\lambda_1|}$ and $u^{\ast}:t\mapsto \cos((\lambda_2-\lambda_1)t)$.
Then, for every $n$ in $\mathbf{N}$, there exists $T_n^{\ast}$ in $(nT^{\ast}-T, nT^{\ast}+T)$ such that
\begin{equation}\label{EQ_estimation_convergence_cosinus}
\frac{1-|\langle \phi_2,\Upsilon^{\frac{u^{\ast}}{n}}(T^{\ast}_n,0) \phi_1 \rangle |}{ 1+ 2 K\| B\|}\leq \frac{ (1+C')\|\pi_2 B \| I}{n},
\end{equation}
with $T^{\ast}=  {\pi}/{2}$, $I={4}/{|\lambda_2-\lambda_1|}$, $ K={2}/{|b_{12}|}$ and
 $$ C'=\sup_{(j,k)\in \Lambda'} \left \{ \left |{\sin \left ( \pi\frac{|\lambda_j-\lambda_k|}{|\lambda_2-\lambda_1|} \right )} \right |^{-1} \right \},$$
where $\Lambda'=\{ (j,k)\in \{1,\ldots,N\}^2 \mbox{ such that } \{j,k\}\cap \{1,2\}\neq \emptyset \mbox{ and } |\lambda_j-\lambda_k|\leq 3/2 |\lambda_2-\lambda_1| \mbox{ and }  b_{jk} \neq 0\}$.
\end{proposition}

\textbf{Proof:}
Fix $\varepsilon>0$, and define $\eta=\varepsilon/K$. We apply the procedure of Section \ref{SEC_Galerkin_approxim} and find $N$ in $\mathbf{N}$ such that $\|(1-\pi_N) B \pi_2\|<\eta$.
By definition of $N$, for every $n$ in $\mathbf{N}$, for every $t\leq T^{\ast}_n$,
$
\|\pi_2 \Upsilon^u_t\phi_1-\pi_2X^{\frac{u^{\ast}}{n}}_{(N)}(t,0)\phi_1\|\leq \varepsilon.
$

A direct computation shows that, for every $\omega>3\pi/T$,
$$
\left |\frac{\int_0^T e^{\mathrm{i}\omega t}u^{\ast}(t)\mathrm{d}t}{\sin(\omega T/2)} \right |=\left |\frac{\omega T^2}{\omega^2T^2-4 \pi^2}\right |\leq\frac{3}{5} \frac{T}{\pi}\leq  I.
$$
 Hence, for every $j,k$ such that $\lambda_{j}-\lambda_{k} >3/2 |\lambda_2-\lambda_1|$, (\ref{EQ_inegalite_non_resonant}) reads
$$\left | \int_0^t h^n_{jk}(\tau)\mathrm{d}\tau \right |\leq \frac{2 |b_{jk}| I}{n}.$$
Estimate (\ref{EQ_estimation_convergence}) becomes
$$
\frac{1-|\langle \phi_2,X^{\frac{u^\ast}{n}}_{(N)}(T^{\ast}_n,0) \phi_1 \rangle |}{ 1+ 2 K\| B\| }\leq
 \frac{ (1+C')\|\pi_2 B \| I}{n},
$$
or, by definition of $N$,
$$
\frac{1-|\langle \phi_2, \Upsilon^{\frac{u^{\ast}}{n}}(T_n^{\ast},0)\phi_1 \rangle|}{ 1+ 2 K\| B\| }\leq 2\varepsilon+ \frac{ (1+C')\|\pi_2 B \| I}{n}.
$$
Conclusion follows when $\varepsilon$ tends to zero.
\qed

Proposition \ref{PRO_cosinus} is the translation in mathematical terms of the well-known fact that the difficulty of inducing a given transition between two eigenstates $j$ and $k$ of the free Hamiltonian is mainly due to eigenstates with energy close to $\lambda_j$ or $\lambda_k$. In other words, in the case of bounded coupling, the convergence in Theorem \ref{PRO_main_result} occurs independently of the  high energy levels of $A$. These levels are rarely known precisely in practice.

\section{Example: Morse potential with truncated dipolar interaction}\label{SEC_examples}

The approximation procedure of Section \ref{SEC_general_case} has been successfully applied to some classical  models, see for instance Boussa\"{i}d, Caponigro \& Chambrion   \cite{FEPS_ACC} and \cite{QG_ACC}  for  the rotation of a 2D molecule, the infinite square potential well (with bounded coupling operator $B$) and a perturbation of the quantum harmonic oscillator (with unbounded potential $B$) .

We present below another type of quantum oscillator. Its main feature is that 
the spectrum of the free Hamiltonian presents a continuous part and a discrete part.Very few controllability results are known for such systems with mixed spectrum (see Mirrahimi \cite{mirrahimi-continuous} for controllability results based on a Lyapunov approach).

\subsection{Modeling}
We consider a diatomic molecule submitted to a time variable electric field with support in some bounded domain. The potential energy of the molecule is modeled with the  Morse potential (see Morse \cite{Morse}). With a suitable choice of units, the Schr\"odinger equation reads
$$
\mathrm{i}\frac{\partial \psi}{\partial t}= -\frac{\partial^2 \psi}{\partial x^2} + V(x) \psi +u(t)  W_M(x) \psi, \quad x>0,
$$
where $V:x\mapsto \alpha^2 \left (e^{-2(x-x_e)}-2e^{-(x-x_e)}\right )$, and $W_M:x\mapsto x$ if $x\leq M$ and $W_M(x)=0$ if $x\geq M$. The positive constants $\alpha$ and $x_e$ describe the physical properties of the molecule: the potential $V$ reaches its minimum at $x_e$, which can be considered as the equilibrium length of the molecule and  $\alpha$ is related to the depth of the potential well. The constant  $M$ (large with respect to $x_e$) is the size of the region where the electric field is assumed to be active.    With the notations of Section \ref{SEC_notations}, $H=L^2((0,\infty),\mathbf{C})$, $A=\mathrm{i}(\Delta-V)$ and $B=-\mathrm{i}W_M$.

Define the integer $N=\lfloor{\alpha-1/2}\rfloor$.  In the following, we assume that $N\geq 4$. The skew-adjoint operator $A$ admits a finite family simple eigenvalues $(-\mathrm{i}\lambda_n)_{0 \leq n \leq N}$, associated with the eigenfuctions $\phi_n:x\mapsto e^{-x} x^{\alpha} P_n(x)$ where $P_n$ is a polynomial function. For every $n\leq N$, $\lambda_n=-(\alpha-n-1/2)^2$. The spectrum of $A$ contains also a continuous part $-\mathrm{i}[0,+\infty)$.

\subsection{Galerkin approximation}
We aim to transfer the system from the first energy level to the second one.
Since the family $(\phi_n)_{n\leq N}$ is not an Hilbert basis of $H$, $(A,B,\mathbf{R})$ does not satisfy Assumption \ref{ASS_assumption_weak} and we cannot apply directly Theorem \ref{PRO_main_result}.

Note that $\langle \phi_{0}, B \phi_{1}\rangle$ tends to $\int_{\mathbf{R}^+} x\phi_0(x) \phi_1(x) \mathrm{d}x \neq 0$ as $M$ tends to infinity. Hence $\langle \phi_{0}, B \phi_{1}\rangle \neq 0$ for $M$ large enough. Moreover, for every $l,m$, $\lambda_m-\lambda_l=(n-m)(-2\alpha+1+(n+m))$.

From now on, we do the generic hypothesis that $\alpha \notin \mathbf{Q}$. In this case,  $\lambda_m-\lambda_l=\lambda_{l'}-\lambda_{m'}$ implies $n-m=n'-m'$ and $n+m=n'+m'$, that is $\{n,m\}=\{n',m'\}$. Moreover, for every $\{j,k\}$,  $|\lambda_j-\lambda_k|\in (\mathbf{N}\setminus\{0\}) (\lambda_1-\lambda_0)$ implies $\{j,k\}=\{0,1\}$.

Fix $\varepsilon>0$. Inspired by Proposition \ref{PRO_cosinus}, we define $u^{\ast}:t\mapsto \cos((\lambda_1-\lambda_0)t)$, $T^{\ast}=\pi/2$, $I=4/|\lambda_1-\lambda_0|$, $K=2/|\langle \phi_0,B\phi_1\rangle|$, $\Lambda'=\{(j,k) \in \{0,1,2,3\}^2| \{j,k\}\neq \{0,1\}\}$ and $$C''=\sup_{(j,k) \in \Lambda'} \left | \sin \left (\pi\frac{\lambda_j-\lambda_k}{\lambda_1-\lambda_0}\right ) \right|^{-1}.$$

From Theorem 2.1, page 525, of \cite{Kato}, for every $\eta>0$, there exists a skew-adjoint operator $A_{\eta}$ such that $A_{\eta}$ admits a complete family of eigenvectors, $A \phi_{l}=A_{\eta}\phi_{l}$  for every $l \leq N$ and $\|A-A_{\eta}\|<\eta$. The pure point spectrum $(-\mathrm{i}\lambda_n^{\eta})_{n\in \mathbf{N}}$ is everywhere dense in $-\mathrm{i}[0,+\infty)$.

The system $(A_{\eta},B,\mathbf{R})$ satisfies Assumption \ref{ASS_assumption_weak}.
For every locally integrable $u$, we denote with $\Upsilon^{u}_\eta$ the propagator of $\frac{d}{dt}\psi=(A_{\eta}+uB_M)\psi$.

We choose an integer $n\geq {(1+C'')M I (1+2 MK)}/{\varepsilon} $
and we define $\eta={\varepsilon}/\left (nT^{\ast}+ T\right ).$  The transition $(0,1)$  of $(A_{\eta},B)$  is non degenerate,  $\lambda_1-\lambda_0=-2(\alpha+1)$ and, for every $l \geq 4$, 
\begin{eqnarray*}
\lambda_l-\left (\lambda_0+\frac{3}{2}(\lambda_1-\lambda_0) \right )&\geq& \lambda_4-\left (\lambda_0+\frac{3}{2}(\lambda_1-\lambda_0) \right )\\
 &\geq &5\alpha-17>0.
\end{eqnarray*}
Proposition \ref{PRO_cosinus} applied to system $(A_{\eta},B,\mathbf{R})$ gives the existence of $T_n^{\ast} \leq nT^{\ast}+T$ such that
$
1-|\langle \phi_1,\Upsilon^{u}_\eta (T_n^{\ast})\phi_0 \rangle| \leq \varepsilon.
$
Since $\|A-A_{\eta}\|<\eta$, we have, for every $t\leq nT^{\ast}+T$,
$$
\left \|\Upsilon^{\frac{u^{\ast}}{n}}(t,0)-\Upsilon^{\frac{u^{\ast}}{n}}_{\eta}(t,0) \right \|<(nT^{\ast}+T)\eta \leq \varepsilon,
$$
and finally
$
1-|\langle \phi_{1}, \Upsilon^{\frac{u^{\ast}}{n}}(T_n^{\ast},0))\phi_{0} \rangle |<2 \varepsilon.
$
\section{Conclusion}
 The contribution of this paper has two elements. First, it proves the validity of the Rotational Wave Approximation for infinite dimensional quantum systems with a pure point spectrum. The result can be partially extended to systems with a mixed spectrum (including a discrete and a continuous part).     The convergence results are accompanied by explicit error estimates which provide explicit time estimates for the approximate controllability of bilinear infinite dimensional quantum systems. Second, a notion of efficiency has been introduced. The efficiency is a measure of the $L^1$-norm of the control needed to achieve a given transition between two eigenstates of the free Hamiltonian. 

The $L^1$-norm of the control appears to play a central role in the analysis of  bilinear quantum systems. 
This is slighty surprising, since one would rather expect that the $L^2$-norm of the control field represents the energy given to the system. It may be an artifact of the semi-classical model that disappears with a more realistic model that takes into account a quantized field.    

 Among other topics, future analysis may concentrate on the design of time-efficient controls or on a systematic treatment  of quantum systems when the free Hamiltonian has a continuous spectrum.
\section*{Acknowledgments}
This research has been supported by the European Research Council, ERC StG 2009 ``GeCoMethods'', contract number 239748, by the ANR project GCM, program ``Blanche'', project number
NT09-504490, and by the Inria Nancy-Grand Est ``COLOR'' project.

It is a pleasure for the author to thank Pierre Rouchon for suggesting him this problem, Ugo Boscain, Marco Caponigro, Julien Salomon, Mario Sigalotti and Dominique Sugny for their valuable input, Nabile Boussa\"{i}d for his most valuable help on functional analysis issues and Denzil Millichap for many corrections.

\bibliographystyle{alpha}
\bibliography{biblioteca}

\end{document}